\documentclass[a4paper,12pt]{article}
\usepackage{amsmath, amsthm, amsfonts, amsxtra, amssymb, psfrag, graphicx, color}
\usepackage{hyperref}

\newtheorem{definition}{Definition}
\newtheorem{theorem}[definition]{Theorem}
\newtheorem{conjecture}{Conjecture}
\newtheorem{proposition}[definition]{Proposition}

\newtheorem{corollary}[definition]{Corollary}
\newtheorem{lemma}[definition]{Lemma}
\newtheorem{fact}[definition]{Fact}
\newtheorem*{claim}{Claim}
\newtheorem{open}{Question}
\graphicspath{{figures/}}
\input{psfig.sty}
\gdef\path{./figures/}

\def\proof{\noindent{\sl Proof.}\ }
\def\qed{\hfill\fbox{\hbox{}}\medskip}
\def\qedclaim{\hfill$\triangle$\smallskip}

\def\P{\mathcal{P}}
\def\A{\mathcal{A}}
\def\B{\mathcal{B}}
\def\C{\mathcal{C}}
\def\D{\mathcal{D}}
\def\F{\mathcal{F}}
\def\mc{\mathcal}
\def\Max{\textmd{Max}}
\def\Min{\textmd{Min}}
\def\AD{A^{\downarrow}}
\def\BD{B^{\downarrow}}
\def\MA{\textmd{Max}(\A)}
\def\da{\downarrow}
\def\tA{\tilde{A}^{\da}}
\def\txt{\textmd}

\title{Linear Extension Diameter of Downset Lattices of 2-Dimensional Posets}
\author{Stefan Felsner and Mareike Massow \\ 
Technische Universit\"at Berlin \\ 
{\tt \{felsner,massow\}@math.tu-berlin.de}}
\date{}

\begin{document}

\maketitle

\begin{abstract}
	The linear extension diameter of a finite poset $\P$ is the maximum distance between a pair of linear extensions of $\P$, where the distance between two linear extensions is the number of pairs of elements of $\P$ appearing in different orders in the two linear extensions. We prove a formula for the linear extension diameter of the Boolean Lattice and characterize the diametral pairs of linear extensions. For the more general case of a downset lattice $\D_\P$ of a 2-dimensional poset~$\P$, we characterize the diametral pairs of linear extensions of $\D_\P$ and show how to compute the linear extension diameter of $\D_\P$ in time polynomial in~$|\P|$.  
\end{abstract}

\section{Introduction}\label{intro}

With a finite poset $\P$ consider its \emph{linear extension graph} $G(\P)$, which has the linear extensions of $\P$ as vertices, with two of them being adjacent exactly if they differ only in a single adjacent transposition. The linear extension graph was originally defined in~ \cite{pruesse91}. It is implicitly used in investigations around finite posets and sorting problems, see e.g.~\cite{bubley-dyer}. 
Explicit research has been initiated in~\cite{reuter}, see also~\cite{naatz}.

The \emph{linear extension diameter} of $\P$, denoted by $\textmd{led}(\P)$, is the diameter of $G(\P)$, see~\cite{felsner-led}. It equals the maximum number of pairs of elements of $\P$ that can be in different orders in two linear extensions of $\P$. A \emph{diametral pair of linear extensions} of $\P$ is a diametral pair of $G(\P)$. 

A realizer of $\P$ is a set $\mathcal{R}$ of linear extensions of $\P$ such that the comparabilities of $\P$ are exactly the intersection of the comparabilities of the linear extensions in $\mathcal{R}$ (cf.~\cite{trotter92}). The dimension of a poset $\P$ is the minimum size of a realizer. If $\P$ is 2-dimensional, i.e., if it has a realizer $\mathcal{R}=\{L_1, L_2\}$, then every incomparable pair $x||y$ of elements of $\P$ appears in different orders in $L_1$ and $L_2$. It follows that $\P$ is 2-dimensional exactly if $\textmd{led}(\P)$ equals the number of incomparable pairs of~$\P$. Figure~\ref{chev} shows a six-element poset $\P$ (called the \emph{chevron}) with its linear extension graph. Note that $\P$ has seven incomparable pairs, but the diameter of its linear extension graph is only six. Hence, the dimension of $\P$ must be at least three. 

   \calc_figscale{50}
    \begin{figure}[htb]
    \centerline{\input{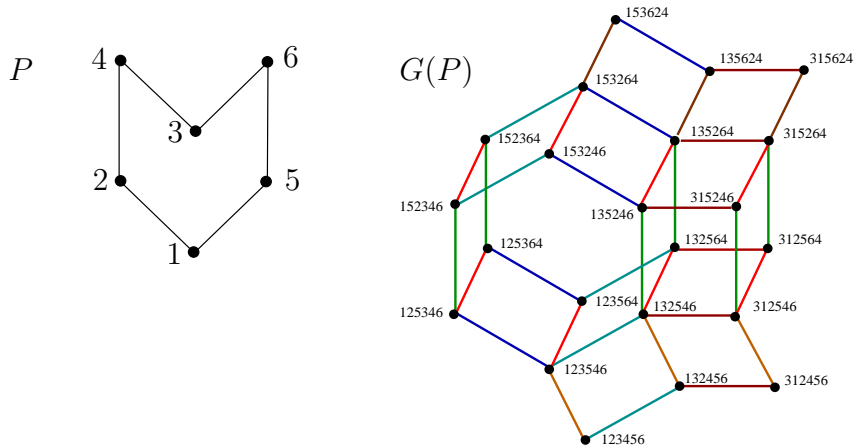}}
    \caption{The chevron and its linear extension graph.\label{fig:chevron-bunt+}}
    \end{figure}
    
\label{chev}

A diametral pair $L_1, L_2$ can be used to obtain an optimal drawing of $\P$: Use $L_1$ and $L_2$ on the two coordinate axes to get a position in the plane for each element of $\P$. Since the number of incomparable pairs which appear in different orders in $L_1$ and $L_2$ is maximized, the resulting drawing has a minimal number of pairs $x,y$ which are comparable in the dominance order, i.e., $x_i < y_i$ for all coordinates $i$, but incomparable in $\P$. See  Figure~\ref{fig_b4} for an example.

 In~\cite{massow-diam} it was shown that, given a poset $\P$, it is NP-complete in general to determine the linear extension diameter of $\P$. In Section~\ref{sec_Bn} we prove a formula for the linear extension diameter of the Boolean lattice, and characterize the diametral pairs of linear extensions of the Boolean lattice. In Section~\ref{sec_latdown} we characterize the diametral pairs of linear extensions for the more general class of downset lattices of 2-dimensional posets. We also show how to compute the linear extension diameter of the downset lattice of a 2-dimensional poset~$\P$ in time polynomial in~$|\P|$. 

Note that the results of Section~\ref{sec_Bn} are contained in the results of Section~\ref{sec_latdown}. 
Section~\ref{sec_Bn} provides an easier acess to our proof techniques. However, since Section~\ref{sec_latdown} is largely self-contained, the self-confident reader may proceed there right away.

\section{Boolean Lattices}
\label{sec_Bn}

Let $B_n$ denote the $n$-dimensional Boolean lattice, that is, the poset on all subsets of $[n]$, ordered by inclusion. In this section we prove the following conjecture from~\cite{felsner-led}.  

\begin{conjecture}[Felsner, Reuter '99]
\label{conj_led_Bn}
	$\textmd{led}(B_n) = 2^{2n-2} - (n+1) \cdot 2^{n-2}$.
\end{conjecture}

We will also characterize the diametral pairs of linear extensions of $B_n$. 
Central to our investigations is a generalization of the reverse lexicographic order. Before defining it, let us clarify some notation: Let $\sigma$ be a permutation of $[n]$. We say that~$i$ is $\sigma$-smaller than~$j$, and write $i <_\sigma j$, if $\sigma^{-1}(i) < \sigma^{-1}(j)$. We define $\sigma$-larger analogously. The $\sigma$-maximum ($\sigma$-minimum) of a finite set~$S$ is the element which is $\sigma$-largest ($\sigma$-smallest) in $S$. For example, if $\sigma=2413$, then $4 <_\sigma 1$, and $\max_\sigma \{3,4\}=3$. 

\begin{definition}
\label{def_sigma_revlex}
Define a relation on the subsets of $[n]$ by setting 
\[
	S<_\sigma T \; \Longleftrightarrow \; \max_\sigma (S \triangle T) \in T
\]
for a pair $S,T \subseteq [n]$. 
\noindent We call this relation the \emph{$\sigma$-revlex order}. 
\end{definition}

In the lemma below we prove that the $\sigma$-revlex order defines a linear extension of $B_n$. We denote this linear extension by~$L_\sigma$. By $\overline{\sigma}$ we denote the reverse of a permutation $\sigma$. In Theorem~\ref{thm_Bn_uniqueness} we prove that the pairs $L_\sigma, L_{\overline{\sigma}}$ are exactly the diametral pairs of linear extensions of $B_n$. Figure~\ref{fig_b4} shows the drawing of $B_4$ resulting from~$L_{\textmd{id}}, L_{\overline{\textmd{id}}}$. 

   \calc_figscale{33}
    \begin{figure}[htb]
    \centerline{\input{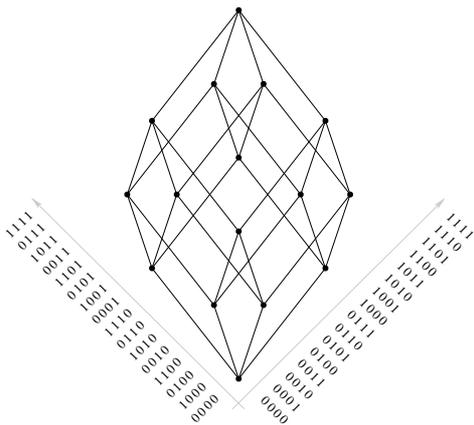}}
    \caption{The drawing of $B_4$ based on the diametral pair $L_{\textmd{id}}, L_{\overline{\textmd{id}}}$.\label{fig:b4-neu}}
    \end{figure}
    
\label{fig_b4}

\begin{lemma}
\label{revlex_LE}
	The relation of being in $\sigma$-revlex order defines a linear extension~$L_\sigma$ of the Boolean lattice. 
\end{lemma}

\proof
We need to show that the $\sigma$-revlex relation defines a linear order on the subsets of $[n]$ which respects the inclusion order. The last part is easy to see: For $S,T \subseteq [n]$ with $S \subseteq T$ we have $S \triangle T \subseteq T$, thus $S<_\sigma T$. 

It is also clear that the relation is antisymmetric and total. So we only need to prove transitivity. Assume for contradiction that there are three sets $A,B,C \subseteq [n]$ with $A <_\sigma B$ and $B <_\sigma C$ and $C <_\sigma A$. Then \mbox{$\max_\sigma (A \triangle B) = b \in B$}, $\max_\sigma (B \triangle C) = c \in C$ and $\max_\sigma (C \triangle A) = a \in A$. 

Note that from $b \in B$ and $c \notin B$ we have that $b \neq c$. Because the situation of the three elements $a,b,c$ is symmetric, it follows that they are pairwise different. Assume that $a = \min_\sigma \{a,b,c\}$. Then by definition of $a$ we know that every element in $C$ which is $\sigma$-larger than $a$ is also in $A$. Hence $c \in A$. Since $c \notin B$, it follows that $c \in A \triangle B$, and thus $\max_\sigma (A \triangle B) = b>_\sigma c$. Now by definition of $c$, every element in $B$ which is $\sigma$-larger than $c$ is also in $C$. Thus $b \in C$. But we also know that $b \notin A$, and hence $b \in C \triangle A$. But this is a contradiction since $b >_\sigma a = \max_\sigma (C \triangle A)$. \qed

\begin{definition}
\label{def_CDI}
For $D \subseteq [n]$ and $I \subseteq [n]\setminus D$ we define~$\C_{D,I}$ as the set of all ordered pairs $(S,T)$ of subsets of $[n]$ with $S \triangle T = D$ and $S \cap T = I$. 
\end{definition}

Clearly, the sets $\C_{D,I}$ partition the ordered pairs of subsets of $[n]$ into equivalence classes. 
Note that these equivalence classes are in bijection with the intervals of $B_n$: The class $\mc{C}_{D,I}$ corresponds to the interval $[I, D \cup I]$, and it contains all pairs $(S,T)$ with $S \cap T = I$ and $S \cup T = D \cup I$.

Another important observation is that we can associate with a subset $X \subseteq D$ the pair \mbox{$(X \cup I, X^c \cup I) \in \mathcal{C}_{D,I}$}, where $X^c = D \setminus X$. On the other hand, for each pair $(S,T) \in \mathcal{C}_{D,I}$ we have $S \setminus I =: X \subseteq D$ and $T = X^c \cup I$ . We obtain the following useful lemma:

\begin{lemma}
	The pairs of a class $\mathcal{C}_{D,I}$ are in bijection with the subsets of $D$. Each class contains  $2^d$ ordered pairs, where $|D|=d$.
\end{lemma}

Recall that the distance between two linear extensions of a poset $\P$ is the number of pairs of elements of $\P$ appearing in different orders in the two linear extensions. We call such a pair a \emph{reversal}. Note that reversals are unordered pairs, whereas the equivalence classes $\C_{D,I}$ contain ordered pairs. Therefore we will often have to switch between ordered and unordered pairs in our proofs. 

The following proposition, settling the lower bound of Conjecture~\ref{conj_led_Bn}, was proved inductively in~\cite{felsner-led}. Here, we give a more combinatorial proof. 

\begin{proposition}
\label{dist_rl_ral}
Given a permutation $\sigma$ of $[n]$, the distance between $L_\sigma$ and $L_{\overline{\sigma}}$ as linear extensions of $B_n$ is
\[
2^{2n-2} - (n+1) \cdot 2^{n-2}.
\]
\end{proposition}
\proof 
We need to count the number of unordered pairs of subsets of $[n]$ appearing in different orders in  $L_\sigma$ and $L_{\overline{\sigma}}$. We claim that an equivalence class $\mc{C}_{D,I}$ contributes exactly $2^{d-2}$ reversals between 
$L_\sigma$ and $L_{\overline{\sigma}}$ if $d \geq 2$, and none if~$d<2$.

Observe that since $\overline{\sigma}$ is the reverse permutation of $\sigma$, the $\sigma$-minimum of a set equals its $\overline{\sigma}$-maximum, and vice versa. Thus we have 
\[
	S <_{\overline{\sigma}} T  \; \Longleftrightarrow \; \min_\sigma (S \triangle T) \in T.
\]

Let $\mathcal{C}_{D,I}$ be an equvialence class as in Definition~\ref{def_CDI}. If $D$ is empty, then~$\C_{D,I}$ contains only the pair $(I,I)$, and thus cannot contribute any reversal. If $D$ consists of only one element, say, $x$, then $\C_{D,I}$ consists of the two pairs $(I, x \cup I)$ and $(x \cup I, I)$. Since $I \subset x\cup I$, the class $\C_{D,I}$ again cannot contribute any reversal. 

So let us assume that $D$ contains at least two elements, and hence $\max_\sigma D \neq \min_\sigma D$. Then a pair in $\mathcal{C}_{D,I}$  corresponding to some $X \subseteq D$ is a reversal between $L_\sigma$ and $L_{\overline{\sigma}}$   if and only if exactly one of the elements $\min_\sigma D$ and $\max_\sigma D$ is contained in $X$. Since we want to count $(S,T)$ and $(T,S)$ only once, let us count the sets $X \subseteq D$ with $\min_\sigma D \in X$ and $\max_\sigma D \notin X$. There are $2^{d-2}$ such sets, and thus $\mc{C}_{D,I}$ contributes $2^{d-2}$ reversals between $L_\sigma$ and $L_{\overline{\sigma}}$ as claimed. 

How many reversals does this yield in total? Each set $D \subseteq [n]$ forms a class together with each set $I \subseteq [n] \setminus D$. Since there are $2^{n-d}$ choices for~$I$, each $D$ with $d \geq 2$ accounts for $2^{n-d} \cdot 2^{d-2}=2^{n-2}$ reversals. There are $2^n-(n+1)$ possibilities to choose a set $D \subseteq [n]$ with $d \geq 2$. Thus the distance between $L_\sigma$ and $L_{\overline{\sigma}}$ is 
\[
2^{n-2}\cdot (2^n-(n+1)) = 2^{2n-2} - (n+1) \cdot 2^{n-2}.
\]
\qed

For proving the other half of the conjecture, we will need Kleitman's Lemma (see~\cite{kleitman-lemma} or~\cite{anderson}):

\begin{lemma}[Kleitman's Lemma]
\label{kleitman}
	Let $\mathcal{F}_1$ and $\mathcal{F}_2$ be families of subsets of $[d]$ which are closed downwards, that is, for every set in $\mc{F}_i$ all its subsets are also in $\mc{F}_i$. Then the following formula holds:
\[
	|\mathcal{F}_1|\cdot|\mathcal{F}_2| \leq 2^d|\mathcal{F}_1 \cap \mathcal{F}_2|.
\]
\end{lemma}


\begin{theorem}
\label{thm_led_Bn}
$\textmd{led}(B_n) = 2^{2n-2} - (n+1) \cdot 2^{n-2}$.
\end{theorem}
\proof
Proposition~\ref{dist_rl_ral} yields $\textmd{led}(B_n) \geq 2^{2n-2} - (n+1) \cdot 2^{n-2}$, since the distance between any pair of linear extensions is a lower bound for the linear extension diameter. To prove that this formula is also an upper bound, we will again use the equivalence classes from Definition~\ref{def_CDI}. Given two linear extensions $L_1, L_2$ of $B_n$, we will show that each $\mathcal{C}_{D,I}$ can contribute at most $2^{d-2}$ reversals between $L_1$ and $L_2$. 

Let us fix a class $\mathcal{C}_{D,I}$. It will turn out that the set $I$ actually plays no role for our argument, therefore we assume that $I = \emptyset$. The only thing this assumption changes is that for a set $X \subseteq D$, now $(X,X^c)$ itself is a pair of $\mathcal{C}_{D,I}$. The reader is invited to check that the following argument goes through unchanged if each $X$ is replaced by $X \cup I$ and each $X^c$ by $X^c \cup I$. 

We say that $X \subseteq D$ is \emph{down} in a linear extension $L$ if $X < X^c$ in $L$. Let~$\mathcal{F}_1$ be the family of subsets of $D$ which are down in $L_1$, and $\mathcal{F}_2$ the family of subsets of $D$ which are down in $L_2$. A pair $(X,X^c) \in \mathcal{C}_{D,I}$ yields a reversal between $L_1$ and $L_2$ exactly if $X$ is down in one $L_i$, but not in the other. Thus our aim is to find an upper bound on $|\mathcal{F}_1 \triangle \mathcal{F}_2|$.

The following key observation captures the essence of transitive forcing between the different pairs: If $X < X^c$ in $L_i$, and $Y \subseteq D$ is a subset of $X$, then $X^c \subseteq Y^c$, and hence by transitivity $Y < X < X^c < Y^c$ in $L_i$. Thus from $X \in\mathcal{F}_i$ it follows $Y \in \mathcal{F}_i$ for every subset~$Y$ of~$X$. This means that $\mathcal{F}_1$ and $\mathcal{F}_2$ each form a family of subsets of $[d]$ which is closed downwards. Hence we can apply Kleitman's Lemma, which yields $|\mathcal{F}_1|\cdot|\mathcal{F}_2| \leq 2^d |\mathcal{F}_1 \cap \mathcal{F}_2|$. 

We observe that for every $L$ and every set $X \subseteq D$, exactly one of~$X$ and~$X^c$ is down in~$L$. Hence we have $|\mathcal{F}_1|=|\mathcal{F}_2| = 2^{d-1}$. It follows that  $|\mathcal{F}_1 \cap \mathcal{F}_2| \geq 2^{d-2}$. 

Also, if $X$ is down in both $L_1$ and $L_2$, then $X^c$ is down in neither. That is, $X \in \mathcal{F}_1 \cap \mathcal{F}_2 \Longleftrightarrow X^c \in (\mathcal{F}_2 \cup \mathcal{F}_1)^c$, and thus $|\mathcal{F}_1 \cap \mathcal{F}_2| = |(\mathcal{F}_2 \cup \mathcal{F}_1)^c |$. From this we obtain
\[
	|\mathcal{F}_1 \triangle \mathcal{F}_2| =2^n-|\mathcal{F}_1 \cap \mathcal{F}_2|- | (\mathcal{F}_2 \cup \mathcal{F}_1)^c | \leq 2^d-2^{d-2}-2^{d-2} = 2^{d-1}.
\]

In $\mathcal{F}_1 \triangle \mathcal{F}_2$, every reversal is counted twice -- once with the set that is down in $L_1$ and once with the set that is down in $L_2$. Therefore the number of (unordered) reversals that $\mathcal{C}_{D,I}$ can contribute is at most $2^{d-2}$. Now we can use the same calculation as in the proof of Proposition~\ref{dist_rl_ral} to show that the total number of reversals is at most $2^{2n-2} - (n+1) \cdot 2^{n-2}$. 
\qed

In the above two proofs we have shown the following fact, which we state explicitely for later reference:

\begin{fact}
\label{lemma_class_contr_Bn}
	If $L, \overline{L}$ is a diametral pair of linear extensions of $B_n$, then each equivalence class $\C_{D,I}$ with $d \geq 2$ contributes exactly $2^{d-2}$ reversals between~$L$ and $\overline{L}$. 
\end{fact}

We have shown that for every permutation $\sigma$ of $[n]$, the linear extensions $L_\sigma$ and $L_{\overline{\sigma}}$ form a diametral pair of the Boolean lattice. Thus we know $n!/2$ diametral pairs. The following theorem proves that these are in fact the only ones.

\begin{theorem}
\label{thm_Bn_uniqueness}
The diametral pairs of the Boolean lattice are unique up to isomorphism. More precisely, if $L,\overline{L}$ is a diametral pair of linear extensions of $B_n$ and $\sigma$ is the order of the atoms in $L$, then $L=L_\sigma$ and $\overline{L}=L_{\overline{\sigma}}$. 
\end{theorem}

\proof
We will show by induction on $k$ that each set of cardinality $k$ is in $\sigma$-revlex order in $L$ and in $\overline{\sigma}$-revlex order in $\overline{L}$ with all sets of cardinality less or equal to $k$. We will use the fact that every equivalence class $\mathcal{C}_{D,I}$ with~$d \geq 2$ contributes exactly $2^{d-2}$ reversals between the diametral linear extensions $L$ and $\overline{L}$. 

Recall that $\sigma$ denotes the order of the atoms in $L$. For every pair $i <_\sigma j$ of atoms, consider the class $\mathcal{C}_{D,I}$ defined by $D=\{i,j\}$ and $I= \emptyset$. This class needs to contribute $2^{2-2}=1$ reversal. Thus, $i$ and $j$ must appear in reversed order in $\overline{L}$. Hence the permutation defining the order of the atoms in $\overline{L}$ is $\overline{\sigma}$. 

Let $L^{k}$ be the restriction of $L$ to the sets of cardinality at most $k$. Our induction hypothesis is that all pairs of sets in $L^{k-1}$ are in $\sigma$-revlex order, that is, $L^{k-1} = L_\sigma^{k-1}$, and that all pairs of sets in $\overline{L}^{k-1}$ are in $\overline{\sigma}$-revlex order, that is, $\overline{L}^{k-1} = L_{\overline{\sigma}}^{k-1}$. For the induction step, we will first show that each set of size $k$ is in the desired order in $L^k$ and $\overline{L}^{k}$  with all sets of strictly smaller cardinality, then that it is in the desired order with all sets of equal cardinality.

So let $A$ be a set of size $k$ in $B_n$ and let $A'$ be the subset of $A$ which is largest in $L^{k-1}$. Let $B$ be the immediate successor of $A'$ in $L^{k-1}$. 

\begin{claim}
	$A$ needs to sit in the \emph{slot} between $A'$ and $B$ in $L^k$.
\end{claim}

Note that $A' = A \setminus a$ for an atom $a \in [n]$.  By induction, all subsets of~$A$ are in $\sigma$-revlex order in $L$. So we know that if $A''$ is a second subset of cardinality $k-1$, then the element of $A$ that is missing in $A'$ is smaller than the element of $A$ that is missing in $A''$. Therefore $a = \min_\sigma A$.

Now observe that since $A' < B$ in $L^{k-1}$ and $|A'|=k-1$ we have $A' || B$. Again by induction we know that $\max_\sigma (A' \triangle B) =b \in B$. If there were $b' \in B \setminus A'$ with $b' \neq b$, then $A' < B \setminus b' < B$ in $L^{k-1}$, which is a contradiction to the choice of $B$. Therefore $B \setminus A' = \{b\}$. 

Because of $A' || B$, there is an element $a' \in A' \setminus B$. We have $b >_\sigma a'$ and therefore also $b >_\sigma a$. Hence, $b=\max_\sigma (A \triangle B)$ and $a = \min_\sigma (A \triangle B)$. 

Consider the class $\mathcal{C}_{D,I}$ defined by $D = A \triangle B$ and $I = A \cap B$. Note that $|D \cup I| = |A \cup B|=|A \cup \{b\}|=k+1$. Choose a set $X \subseteq D \setminus \{a,b\}$. Then $|X \cup I| \leq k-1$, thus we can apply the induction hypothesis to get $X \cup I < b \, \cup I$ in $L$, and with $b \in X^c$ it follows $X \cup I < X^c \cup I$ in $L$. Analogously we have $X \cup I < a \cup I < X^c \cup I$ in $\overline{L}$. Thus the pair $(X, X^c)$ does not yield a reversal between $L$ and $\overline{L}$, and neither does the pair $(X^c, X)$. 

There are $2^{d-2}$ choices for $X$, hence we have found $2^{d-1}$ ordered pairs in~$\C_{D,I}$ which do not yield reversals between~$L$ and $\overline{L}$. By Lemma~\ref{lemma_class_contr_Bn}, the class~$\C_{D,I}$ contributes exactly $2^{d-2}$ reversals. But the remaining $2^{d-1}$ ordered pairs in $\C_{D,I}$ can yield at most $2^{d-2}$ reversals, thus, they all have to be reversed between $L$ and $\overline{L}$. 
It follows that all subsets $Y$ of $D$ containing exactly one of the two atoms $a$ and $b$ have to be down in exactly one of the two linear extensions. 

In particular, we can choose $Y = A \setminus I$ to see that $\{A,B\}$ must be a reversal. In~$\overline{L}$, we know the order of~$A$ and~$B$: Set \mbox{$A''= I \cup a \subseteq A$}, then \mbox{$\min_\sigma (B \triangle A'') = a \in A''$}. But this means \mbox{ $\max_{\overline{\sigma}} (B \triangle A'') = a \in A''$}. So we have  $B < A'' $ in~$\overline{L}$ by induction and thus $B< A$ in~$\overline{L}$ by transitivity. Hence it follows that $A < B$ in $L$. This proves our claim that $A$ has to sit in the slot between $A'$ and $B$ in $L$. \qedclaim

Because $\max_\sigma(A \triangle B) = b \in B$, by showing $A<B$ in $L^k$ we have shown that $A$ is in $\sigma$-revlex order with $B$ in $L^{k}$. Since the slot after $A'$ is the lowest possible position for $A$ in $L^k$, it follows by transitivity that $A$ is in $\sigma$-revlex order with all sets of cardinality $<k$ in $L^k$. By reversing the roles of $L$ and $\overline{L}$ we see that $A$ also has to be in $\overline{\sigma}$-revlex order in $\overline{L}$ with all sets of smaller cardinality. Now we will show that all pairs of sets with equal cardinality~$k$ need to be in $\sigma$-revlex order in $L$. 

Let us consider two sets $A_i, A_j \in B_n$ with cardinality $k$. If they are inserted into different slots in $L^k$, then their order in $L^k$ equals their order in~$L_{\sigma}$, thus they are in $\sigma$-revlex order. If they go into the same slot, this means that they have the same largest $(k-1)$-subset $A'$ in $L$. Thus $|A_i \triangle A_j|=2$. We also know that $a_i:= A_i \setminus A' = \min_\sigma A_i$ and $a_j:=A_j \setminus A' = \min_\sigma A_j$. Assume without restriction that $a_i<_\sigma a_j$. Then for the pair of sets $\{A_j, \{a_i\}\}$ we know $a_i=\min_\sigma (A_j \triangle \{a_i\})=\max_{\overline{\sigma}} (A_j \triangle \{a_i\})$. Thus by induction, \mbox{$A_j < a_i < A_i$} in~$\overline{L}$, that is, $A_i$ and $A_j$ belong to different slots in $\overline{L}$. Now since the class $\mathcal{C}_{D,I}$ containing $(A_i,  A_j)$ needs to contribute $2^{2-2}=1$ reversal between $L$ and $\overline{L}$, and $(A_i,  A_j)$ is the only incomparable pair in this class (except for the reversed pair), we know that we must have $A_i < A_j$ in $L^k$. Since $\max_\sigma (A_i \triangle A_j) = a_j$, this means that $A_i$ and $A_j$ are in $\sigma$-revlex order in $L^k$.

Again we can apply the same argument with the roles of $L$ and $\overline{L}$ reversed to show that all sets of cardinality $k$ are in $\overline{\sigma}$-revlex order in $\overline{L}$, thus $\overline{L}^{k} = L_{\overline{\sigma}}^{k}$. By induction we conclude that $L=L_\sigma$ and $\overline{L}=L_{\overline{\sigma}}$. \qed

Another conjecture in ~\cite{felsner-led} suggested a connection between diametral pairs and reversed critical pairs. A \emph{critical pair} of a poset $\P$ is an ordered pair~$(x,y)$ of elements of $\P$ such that all elements smaller than $x$ are also smaller than~$y$ in $\P$, and all elements larger than~$y$ are also larger than~$x$ in~$\P$. Critical pairs appear as an important ingredient in the dimension theory of posets, see for example~\cite{trotter92}. A critical pair $(x,y)$ is \emph{reversed} in a linear extension $L$ of $\P$ if $x> y$ in $L$. 

It was conjectured in~\cite{felsner-led} that in every diametral pair of linear extensions of $\P$, at least one of the two linear extensions reverses a critical pair of elements of $\P$. In~\cite{massow-diam} it was shown that the conjecture is false in general, but that almost all posets have the stronger property of being \emph{diametrally reversing}, which means that \emph{every} linear extension contained in a diametral pair reverses a critical pair. Still for Boolean lattices it remained open even whether they have the weaker property. 

With the help of Theorem~\ref{thm_Bn_uniqueness}, we can now settle this question:

\begin{corollary}
	Boolean lattices are diametrally reversing.
\end{corollary}

\proof
Let $L$ be a linear extension of $B_n$ that is contained in a diametral pair. Then by Theorem~\ref{thm_Bn_uniqueness} we know that $L=L_\sigma$ for some permutation~$\sigma$ of~$[n]$. It can be checked that the $n$ atom-coatom pairs $( \{i\}, [n] \setminus i )$ for $i \in [n]$ are critical pairs of $B_n$ (in fact, these are the only ones). Now let $i=\max_\sigma [n]$. Then $\max_\sigma (\{i\} \triangle [n] \setminus i) = i$. Thus $[n] \setminus i < \{i\}$ in $L_\sigma$, which means that $L_\sigma$ reverses this critical pair. Hence every linear extension of $B_n$  contained in a diametral pair reverses a critical pair of $B_n$. \qed

\section{Downset Lattices of 2-Dimensional Posets}
\label{sec_latdown}

The Boolean lattice $B_n$ can be viewed as the distributive lattice of downsets of the $n$-element antichain. Now let $\P$ be an arbitrary poset. Denote by~$\mathcal{D}_\P$ the downset lattice of $\P$, that is, the poset on all downsets of $\P$, ordered by inclusion. In Section~\ref{subsec_latdown_char} we give an upper bound for the linear extension diameter of $\mc{D}_\P$ and show that it is tight if $\P$ is 2-dimensional. We also 
characterize the diametral pairs of linear extensions of $\mathcal{D}_\P$ for 2-dimensional~$\P$. For the proofs, we make use of the main ideas from the last section. In Section~\ref{subsec_latdown_comp} we show how to compute the linear extension diameter of  $\mathcal{D}_\P$ in time polynomial in $|\P|$.  

We can use our results to obtain optimal drawings of $\mathcal{D}_\P$ as described in the introduction, see Figure~\ref{fig_downset_latt} for an example.

   \calc_figscale{33}
    \begin{figure}[htb]
    \centerline{\input{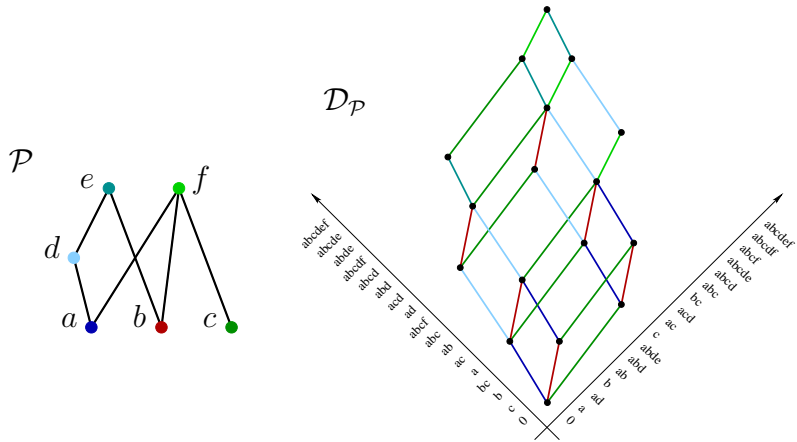}}
    \caption{A drawing of a downset lattice based on a diametral pair.\label{fig:distr-lat-of-2d}}
    \end{figure}
    
\label{fig_downset_latt}

In this section we make fundamental use of the canonical bijection between downsets and antichains of a poset. We  frequently switch back and forth between the two viewpoints, concentrating on the one or the other in our proofs. We write $\AD$ to refer to the downset generated by an antichain~$A$. We write $\MA$ to refer to the antichain of maxima of a downset $\A$.

\subsection{Characterizing Diametral Pairs}
\label{subsec_latdown_char}

This subsection is organized as follows: We first prove an upper bound for the linear extension diameter of arbitrary downset lattices, using a generalization of the equivalence classes from Definition~\ref{def_CDI}. Then we show that for downset lattices of 2-dimensional posets, a generalization of the revlex orders attains this bound. Finally we prove that these are the only pairs of linear extensions attaining the bound. 

\begin{definition}
\label{def_equ_classes_gen}
Let $\P$ be a poset, $D \subset \P$ and $I \subset \P \setminus D$. We define $\C_{D,I}$ as the set of all ordered pairs $(A,B)$ of antichains of $\P$ with $D= A \triangle B$ and~$I = A \cap B$. 
\end{definition}

It is easy to see that the sets $\C_{D,I}$ partition the ordered pairs of antichains of $\P$ into equivalence classes. Note that for a class $\C_{D,I}$, the sets $D$ and $I$ are disjoint. Furthermore, there is no relation in $\P$ between any element of~$I$ and any element of $D$.

\begin{lemma}
\label{lemma_equ_p[d]}
Let $\C_{D,I}$ be an equivalence class as defined above. Let  $\P[D]$ be the subposet of $\P$ induced by the elements of $D$, and let $\mc{K}$ be the set of connected components of $\P[D]$. Then the pairs in $\C_{D,I}$ are in bijection with the subsets of $\mc{K}$. 
\end{lemma}

\begin{figure}[ht]
  \begin{center}
    \psfrag{k1}{$\kappa_1$}
    \psfrag{k2}{$\kappa_2$}
    \psfrag{k3}{$\kappa_3$}
    \psfrag{A}{$\in A$}
    \psfrag{B}{$\in B$}
    \includegraphics[scale=0.5]{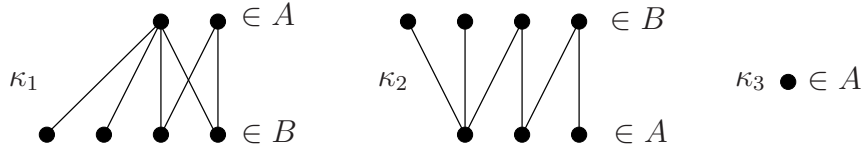}
    \label{fig_PinducedbyD}
\caption{An example for $\P[D]$ with three components. The assignment of minima and maxima to $A$ and $B$ specifies one of the eight pairs in the class.
}
   \end{center}
\end{figure}

\proof
For a given equivalence class $\C_{D,I}$, let $(A, B)$ be a pair in $\C_{D,I}$, thus we have $D = A \triangle B$. First observe that since $A$ and  $B$ are antichains, $\P[D]$ is a poset of height at most 2 (see Figure~\ref{fig_PinducedbyD}). Thus all elements of $D$ belong to the antichain $\Max(D)$ of maxima of $\P[D]$ or to the antichain $\Min(D)$ of minima of $\P[D]$. 

Consequently, every connected component~$\kappa$ of~$\P[D]$ is either a single element or has height $2$. If $\kappa$ is a single element, it belongs either to~$A$ or to~$B$. If  $\kappa$ has height~2, there are again two possibilities: Either the maxima of $\kappa$ all belong to $A$ and its minima all belong to $B$, or the maxima of $\kappa$ all belong to $B$, and its minima all belong to $A$. Note that in the second case the minima also belong to $\BD$, and in the first case also to $\AD$. 

We can get a different pair $(A', B') \in \C_{D,I}$ by switching the roles of $A$ and~$B$ in one component of $\P[D]$. We can do this switch independently for each component, but we have to do it for the whole component to ensure that the elements belonging to $A$ and respectively $B$ still form an antichain. 

Let $K \subseteq \mc{K}$ be a subset of the components of $\P[D]$. With $K$ we associate a subset $X_K$ of $D$ by setting
\[
	X_K = \bigcup_{\kappa \in K} \textmd{Max}(\kappa) \cup \bigcup_{\kappa \notin K, \, |\kappa|>1} \textmd{Min}(\kappa). 
\]
Let $X_K^c = D \setminus X_K$. Now define a map from the powerset of $\mc{K}$ to  $\C_{D,I}$ via
\[
	K \mapsto (A_K, B_K)=(X_K \cup I,  X^c_K \cup I). 
\]
We claim that this defines a bijection. 
Recall that $D$ and $I$ are disjoint and that there are no relations between them. It follows that 
\[
A_K \triangle B_K = (X_K \cup I) \triangle (X^c_K \cup I) = X_K \cup X_K^c = D,
\]
and
\[
A_K \cap B_K =  (X_K \cup I) \cap (X^c_K \cup I)  = I.
\]
Thus $K$ is indeed mapped to a pair in $\C_{D,I}$. 

If $K$ and $K'$ are two different subsets of $\mc{K}$, then $X_K$ and $X_{K'}$ differ on at least one component of $\P[D]$, and hence $(A_K, B_K) \neq (A_{K'}, B_{K'})$. Therefore our map is injective. Also, given a pair $(A,B) \in \C_{D,I}$ it is easy to construct $K \subseteq \mc{K}$ with $(A, B)= (A_K, B_K)$. Thus our map is a bijection as claimed.
\qed

Note that the above lemma implies that there are exactly $2^d$ pairs in the class~$\C_{D,I}$, where $d=|\mc{K}|$ denotes the number of connected components of~$\P[D]$.

For the following, let us keep in mind that $\C_{D,I}$ contains ordered pairs, whereas reversals, constituting the distance between two linear extensions, are unordered pairs. 

\begin{theorem}
\label{thm_diam_downlat_upper}
Let $\mc{D}_\P$ be the downset lattice of a poset $\P$. The linear extension diameter of $\mc{D}_\P$ is bounded by a quarter of the number of pairs $(A, B)$ of antichains of $\P$ such that $\P[A \triangle B]$ has at least two connected components.
\end{theorem}

\proof 
Let $L_1, L_2$ be an arbitrary pair of linear extensions of $\mathcal{D}_\P$ and $\mathcal{C}_{D,I}$ an equivalence class as in Definition~\ref{def_equ_classes_gen}.
First note that if $D$ is empty, then the class $\mathcal{C}_{D,I}$ only consists of a single pair, namely, $(I,I)$. This class cannot contribute any reversal. Now let $D$ be non-empty, and assume that $\P[D]$ is connected. Then with Lemma~\ref{lemma_equ_p[d]} we know that $\mathcal{C}_{D,I}$ consists of the two pairs $(A, B)$ and $(B, A)$, where $A = \Max(D) \cup I$ and $B =(D \setminus \Max(D)) \cup I$. But then we have $\BD \subset \AD$, that is, the two downsets form a comparable pair in~$\mc{D}_\P$. Hence~$\mathcal{C}_{D,I}$ cannot contribute any reversals if $\P[D]$ is connected.

Now let us assume that $\P[D]$ has at least two components. We claim that at most half of the pairs contained in $\C_{D,I}$ can be reversed between $L_1$ and~$L_2$. This means that the number of (unordered) reversals that each class can contribute is at most a quarter of the number of pairs that it contains. 

 Recall the notions from Lemma~\ref{lemma_equ_p[d]}. For a pair $(A,B) \in \C_{D,I}$, let us call~$\AD$ \emph{down in $L_i$} if $\AD < \BD$ in $L_i$, for $i=1,2$. The pairs in $\mathcal{C}_{D,I}$ are in bijection with the subsets of $\mc{K}$. Define the family $\mathcal{F}_i$ as the family of subsets $K$ of $\mc{K}$ such that $A_K^{\downarrow}$ is down in $L_i$. 

A pair $\{\AD,\BD\}$ is a reversal between $L_1$ and $L_2$ exactly if $\AD$ is down in one of the two linear extensions, but not in both. Put differently, a pair $(A_K,B_K)$ yields a reversal exactly if $K$ is contained in one of the~$\mc{F}_i$, but not in both. Thus we are interested in bounding $|\mathcal{F}_1 \triangle \mathcal{F}_2|$.

Let $K \in \F_i$, and $K' \subset K$. We claim that  $K' \in F_i$. Indeed, by definition of $X_K$ we have $(X_{K'} \cup I)^{\downarrow} \subset (X_K \cup I)^{\downarrow }$, that is, $A_{K'}^{\downarrow} \subset A_K^{\downarrow}$. Analogously it holds $(X_{K}^c \cup I)^{\downarrow} \subset (X_{K'}^c \cup I)^{\downarrow }$ and hence $B_K^{\downarrow} \subset B_{K'}^{\downarrow}$. 
Thus from $A_K^{\downarrow} < B_K^{\downarrow}$ in $L_i$ it follows that \mbox{$A_{K'}^{\downarrow} < A_K^{\downarrow} < B_K^{\downarrow} < B_{K'}^{\downarrow}$} \mbox{in $L_i$}. Therefore the families $\mathcal{F}_1$ and $\mathcal{F}_2$ are closed downwards, and we can apply Kleitman's Lemma as in the proof of Theorem~\ref{thm_led_Bn}. This yields \mbox{$|\mathcal{F}_1 \cap \mathcal{F}_2| \geq 2^{-d} \cdot |\F_1|\cdot |\F_2|$}, where $d=|\mc{K}|$. 

Observe that if $(A,B)  \in \C_{D,I}$ is associated with the set $K \subseteq \mc{K}$, then $(B, A)$ is associated with the set $K^c = \mc{K} \setminus K$. Now in each $L_i$, either $\AD$ is down or $\BD$ is down. Thus $K \in \F_i$ exactly if $K^c \notin \F_i$. Hence \mbox{$|\F_1|= |\F_2| = 2^{d-1}$}. It follows that $|\mathcal{F}_1 \cap \mathcal{F}_2| \geq 2^{d-2}$. 
Similarly, if $\AD$ is down in both $L_i$, then~$\BD$ is down in neither, and vice versa. This means that $K \in \F_1 \cap \F_2$ exactly if $K^c  \in (\F_1 \cup \F_2)^c$. Therefore $|\mathcal{F}_1 \cap \mathcal{F}_2| = | (\F_1 \cup \F_2)^c|$. We conclude that 
\[
	|\mathcal{F}_1 \triangle \mathcal{F}_2| = 2^d - |\mathcal{F}_1 \cap \mathcal{F}_2| - | (\F_1 \cup \F_2)^c| \leq 2^d - 2^{d-2} - 2^{d-2} = 2^{d-1}. 
\]

We have thus shown that each class $\C_{D,I}$ can contribute at most $2^{d-2}$ reversals between two linear extensions of $\D_\P$ as claimed. 
\qed

Next we are going to prove that the bound of the above theorem is tight in the case of a 2-dimensional poset $\P$. To do so, we define a particular type of linear extension of $\mc{D}_\P$.  Let $\sigma$ be a linear extension of $\P$. For a set $S \subseteq \P$, let $\max_\sigma S$ be the element of $S$ which is largest in $\sigma$. In analogy to Definition~\ref{def_sigma_revlex}, we have:

\begin{definition}
	We define a relation $<_\sigma$ on the pairs $\{\AD,\BD\}$ of downsets of $\P$ as follows: 
\[
	\AD <_\sigma \BD\; \Longleftrightarrow \; \max_\sigma (A \triangle B)  \in B.
\]
We call this relation the \emph{$\sigma$-revlex order}.
\end{definition}

\begin{lemma}
	Let $\P$ be a poset and $\sigma$ a linear extension of $\P$. The relation of being in $\sigma$-revlex order defines a linear extension~$L_\sigma$ of $\mc{D}_\P$. 
\end{lemma}
\proof
Since $\max_\sigma (\AD \triangle \BD)=\max_\sigma (A \triangle B)$, we can equivalently define the relation $<_\sigma$ directly via the downsets. Thus this relation is just a restriction of the $\sigma$-revlex order on all subsets of $\P$ (see Definition~\ref{def_sigma_revlex}) to the downsets of $\P$. In Lemma~\ref{revlex_LE}, we proved that the $\sigma$-revlex order on all subsets is a linear order and that it extends the inclusion relation. This carries over to the restriction to the downsets of~$\P$. 
\qed

Note that for the special case~$\P=B_n$, the above lemma was already shown in~\cite{hosten-morris}. 

Now let $\sigma$ be a linear extension of $\P$ which is contained in a diametral pair of~$\P$. Since $\P$ is 2-dimensional, $\sigma$ has a unique partner $\overline{\sigma}$ with which it forms a realizer, i.e., a diametral pair. Thus all incomparable pairs of $\P$ are reversals between $\sigma$ and~$\overline{\sigma}$. 
We will frequently use the following characterization, which is implicite already in~\cite{dushnik41}, cf.~also~\cite{moehring-tractable}. 

\begin{lemma}
	Let $\P$ be a poset of dimension $2$. A linear extensions~$\pi$ of~$\P$ is contained in a realizer of cardinality~2 exactly if it is  \emph{non-separating}, that is, for $u,v \in \P$ with $u < v$ in~$\P$, there is no element $x \in \P$ with $x \, || \, u,v$ in~$\P$ and $u < x < w$ in $\pi$. 
\end{lemma}

In particular, $\sigma$ and $\overline{\sigma}$ are non-separating. If a linear extension is not non-separating, we call it \emph{separating}. 

\begin{theorem}
\label{thm_Lsigma_diam}
Let $\P$ be a 2-dimensional poset, and let $\sigma, \overline{\sigma}$ be a diametral pair of linear extensions of $\P$. Then $L_\sigma, L_{\overline{\sigma}}$ is a diametral pair of linear extensions  of $\mathcal{D}_\P$. 
\end{theorem}

\proof We will show that for each class $\mathcal{C}_{D,I}$, the pair $L_\sigma, L_{\overline{\sigma}}$ realizes the maximum possible number of reversals. 
Fix an arbitrary class $\mathcal{C}_{D,I}$. We have seen that $\mathcal{C}_{D,I}$ cannot contribute any reversals if $\P[D]$ is empty or consists of only one component. Now let $\mc{K}$ be the set of connected components of $\P[D]$ as before, and assume that $|\mc{K}| =d \geq 2$. 

A pair $(A, B) \in \C_{D,I}$ yields a reversal between $L_\sigma$ and $ L_{\overline{\sigma}}$ exactly if one of the two elements $\max_\sigma (A \triangle B)$ and $\max_{\overline{\sigma}} (A \triangle B)$ is contained in~$A$, and the other in $B$. Since $(A,B)$ and $(B,A)$ contribute only one reversal, let us count the number of pairs $(A, B)$ with \mbox{$\max_\sigma (A \triangle B) \in A$} and \mbox{$\max_{\overline{\sigma}} (A \triangle B) \in B$}. 

Note that the $\sigma$-maximum of $A \triangle B$ belongs to the antichain $\Max(D)$, which is completely reversed between $\sigma$ and~$\overline{\sigma}$. It follows that we have $\max_\sigma \Max(D) = \min_{\overline{\sigma}} \Max(D)$ and  $\max_{\overline{\sigma}} \Max(D) = \min_\sigma \Max(D)$. We claim that $\max_\sigma \Max(D)$ and $\min_\sigma \Max(D)$ lie in different components of~$\P[D]$. 

Suppose for contradiction that $\max_\sigma \Max(D)$ and $\min_\sigma \Max(D)$ both lie in the component $\kappa$ of $\P[D]$. Since we assumed that $\P[D]$ has at least two components, we can choose an element $x \in \Max(D)$ which is not contained in~$\kappa$ and thus has no relation to any element of~$\kappa$. Then \mbox{$\min_\sigma \Max(D) < x < \max_\sigma \Max(D)$} in $\sigma$. Denote by $\kappa_1$ the set of elements of $\kappa$ which are $\sigma$-smaller than $x$, and by $\kappa_2$ the set of elements of $\kappa$ which are $\sigma$-larger than $x$. Both sets are non-empty since \mbox{$\min_\sigma \Max(D) \in \kappa_1$} and $ \max_\sigma \Max(D) \in \kappa_2$. Since $\kappa$ is a connected component of $\P[D]$, there are elements $u \in \kappa_1$ and $v \in \kappa_2$ with $u < v$ in $\P$. But then $u < x < v$ in $\sigma$ with $x || u,v$ in $\P$, and this is a contradiction because $\sigma$ is non-separating.

We have shown that $\max_\sigma (A \triangle B)$ and $\max_{\overline{\sigma}} (A \triangle B)$ lie in different components of $\P[D]$, say, $\max_\sigma (A\triangle B) \in \kappa$ and $\max_{\overline{\sigma}} (A \triangle B) \in \lambda$ with $\kappa, \lambda \in \mc{K}$. Now let \mbox{$K \subseteq \mc{K}$} be a set of components with $\kappa \in K$ and $\lambda \notin K$. Consider the pair $(A_K, B_K) \in \C_{D,I}$. By definition, we have $\max_\sigma (A \triangle B) \in A_K$ and $\max_{\overline{\sigma}} (A \triangle B) \in B_K$. Thus the pair $(A_K, B_K)$ contributes a reversal between $L_\sigma$ and $ L_{\overline{\sigma}}$. There are $2^{d-2}$ possibilities of choosing $K$. Thus every class~$\C_{D,I}$ with $d \geq 2$ contributes at least $2^{d-2}$ reversals between $L_\sigma$ and $ L_{\overline{\sigma}}$. We have seen in Theorem~\ref{thm_diam_downlat_upper} that this is maximal. Therefore $L_\sigma$ and $ L_{\overline{\sigma}}$ are a diametral pair of linear extensions of~$\mathcal{D}_\P$. 
\qed

We have now proved the fact that for a 2-dimensional poset $\P$, every class~$\C_{D,I}$ for which $\P[D]$ has at least two components contributes exactly $2^{d-2}$ reversals between a diametral pair of linear extensions of $\mc{D}_\P$. Since there are $2^d$ pairs in $\C_{D,I}$, we have the following corollary: 

\begin{corollary}
\label{corr_led_downlat}
Let $\mc{D}_\P$ be the downset lattice of a 2-dimensional poset $\P$. The linear extension diameter of $\mc{D}_\P$ equals a quarter of the number of pairs $(A, B)$ of antichains of $\P$ such that $\P[A \triangle B]$ has at least two connected components.
\end{corollary}

Note that the proof of Conjecture~\ref{conj_led_Bn}, thus, the formula for the linear extension diameter of the Boolean lattice, follows from this corollary. A special case of the proof of Theorem~\ref{thm_led_factorlattice} yields $\txt{led}(B_n)= 2^{2n-2}-(n+1)2^{n-2}$. 


In the next theorem we characterize the diametral pairs of linear extensions of downset lattices of 2-dimensional posets.

\begin{theorem}
\label{thm_unique_downlat}
	Let $\P$ be a 2-dimensional poset, and let $L,\overline{L}$ be a diametral pair of linear extensions of $\mathcal{D}_\P$. Let $\sigma$ be the linear extension of $\P$ defined by the order of the downsets $x^{\downarrow}$ for $x \in \P$ in $L$. Then $\sigma$ is contained in a diametral pair $\sigma, \overline{\sigma}$ of linear extensions of $\P$, and we have $L=L_\sigma$ and $\overline{L}= L_{\overline{\sigma}}$. 
\end{theorem}

\proof
We again use the equivalence classes from Definition~\ref{def_equ_classes_gen}. For each incomparable pair $x,y \in \P$, consider the equivalence class $\mathcal{C}_{D,I}$ defined by  $D=\{x,y\}$ and $I=\emptyset$. Then $\P[D]$ consists of two singletons. Since $L,\overline{L}$ is a diametral pair, this class must contribute $2^{2-2}=1$ reversal, and thus $x^{\downarrow}$ and $y^{\downarrow}$ must appear in opposite order in $L$ and $\overline{L}$. Recall that $\sigma$ is the linear extension of $\P$ defined by the order of the downsets $x^{\downarrow}$ in $L$. Let us denote the linear extension defined analogously for $\overline{L}$ by $\overline{\sigma}$. We have seen that every incomparable pair of elements of $\P$ must be a reversal between $\sigma$ and $\overline{\sigma}$. Hence, $\sigma, \overline{\sigma}$ is a diametral pair of linear extensions of~$\P$. 

In the following, we use induction on the cardinality of the downsets of~$\P$ to show that $L=L_\sigma$ and $\overline{L}=L_{\overline{\sigma}}$, in analogy to the proof of Theorem~\ref{thm_Bn_uniqueness}. More precisely, we show that each downset of cardinality $k$ is in $\sigma$-revlex order in $L$ and in $\overline{\sigma}$-revlex order in $\overline{L}$ with all downsets of cardinality $\leq k$, by induction on $k$. We use the fact that every equivalence class $\mathcal{C}_{D,I}$ for which $\P[D]$ is disconnected contributes exactly $2^{d-2}$ reversals between $L$ and $\overline{L}$. Note that we have settled the base case already: All downsets of cardinality~$1$ are of the form $x^{\downarrow}$ for some minimal element $x \in \P$, and we have shown in the previous paragraph that these behave as expected.

Let $L^{k}$ be the restriction of $L$ to the sets of cardinality at most $k$. Our induction hypothesis is that $L^{k-1} = L_\sigma^{k-1}$ and $\overline{L}^{k-1} = L_{\overline{\sigma}}^{k-1}$. We structure the induction step as follows: We first show that each set of size $k$ is in the desired order in $L^k$ and $\overline{L}^k$ with all sets of smaller size. This will be the main part of the proof. Then we show that all pairs of sets of equal size~$k$ are in the desired order  in $L^k$ and $\overline{L}^k$.

So let $\AD$ be a downset of cardinality $k$ of $\P$. Let $\tA$ be the subset of $\AD$ which is largest in $L^{k-1}$, and let $\BD$ be its immediate successor in $L^{k-1}$. 

\begin{claim}
	$\AD$ needs to sit in the \emph{slot} between $\tA$ and $\BD$ in $L^k$.
\end{claim}

Proving this claim requires some technical details. Here is an outline of what we are going to do: We first locate the elements $\max_\sigma (\AD \triangle \BD)=\max_\sigma (A \triangle B)$ and \mbox{$ \max_{\overline{\sigma}} (\AD \triangle \BD)=\max_{\overline{\sigma}} (A \triangle B)$}. Using these we see that $\{\AD,\BD\}$ needs to be a reversal between $L$ and $\overline{L}$. From the order of $\AD$ and~$\BD$ in $\overline{L}$ we can finally deduce that $\AD < \BD$ in~$L^k$.

We have $\tA = \AD \setminus a$ for some $a \in A$. All subsets of $\AD$ are in \mbox{$\sigma$-revlex} order in $L$ by induction. So we know that if $\hat{A}^{\da}$ is a second subset of cardinality~$k-1$ of $\AD$, then the element of $\AD$ that is missing in $\tA$ is $\sigma$-smaller than the element of $\AD$ that is missing in $\hat{A}^{\da}$. Thus we can conclude that $a = \min_\sigma A$. Since the antichain $A$ is completely reversed between $\sigma$ and $\overline{\sigma}$, it follows that \mbox{$a = \max_{\overline{\sigma}} A$}. 

Now observe that since $\tA < \BD$ in $L^{k-1}$ and $|\tA|=k-1$ we have $\tA \, || \, \BD$. By induction we know that $\max_\sigma (\tA \triangle \BD) \in \BD$. Let $b$ be the $\sigma$-smallest element of $\BD \setminus \tA$ which is $\sigma$-larger than all elements of $\tA \setminus \BD$. Then $(\tA \cap \BD) \cup b$ is a downset of $\P$. Observe that by induction $\tA < (\tA \cap \BD) \cup b$ in~$L^{k-1}$. Since $(\tA \cap \BD) \cup b \subseteq \BD$ we must have $(\tA \cap \BD) \cup b=\BD$ by the choice of $\BD$. Thus $\BD \setminus \tA=\{b\}$ and $\max_\sigma (\tA \triangle \BD)= b$. We will use the next three paragraphs to show that $ \max_\sigma (\AD \triangle \BD)=b$ and $ \max_{\overline{\sigma}} (\AD \triangle \BD)=a$.

Note that $\BD \not\subseteq \AD$ by the choice of $\tA$, and thus $\BD \setminus \AD = \{b\}$. Hence $b \not\leq a$ in $\P$. On the other hand, $a \not< b$ in $\P$ because otherwise $a \in \BD \setminus \tA$ and thus $a = b$, a contradiction. It follows that $a||b$ in~$\P$. 

Next let us show that $a <_\sigma b$.
Because of $|\tA|=k-1$, we know that $\tA \setminus \BD \neq \emptyset$. Let $a' \in \tilde{A} \setminus B$. If $a || a'$ in $\P$, then $a' \in A \setminus B$ and with $a= \min_\sigma A$ it follows that $a <_\sigma a'$. Together with $a' <_\sigma b$ we get $a <_\sigma b$. 
If $a$ and $a'$ are comparable, then $a'< a$ in~$\P$. Now assume for contradiction that $b <_\sigma a$. Then we have $a' <_\sigma b <_\sigma a$, with $b \, || \, a', a$ and $a'< a$ in~$\P$. This means that $\sigma$ is a separating linear extension. But since $\sigma$ is contained in the realizer $\sigma, \overline{\sigma}$ of $\P$, this is a contradiction. 

We have shown that $a <_\sigma b$. We knew already that  \mbox{$\max_\sigma (\tA \triangle \BD)= b$}, and because $a$ is the only element in $\AD \setminus \tA$, we can conclude that \mbox{$\max_\sigma (\AD \triangle \BD)=b$}. 
Also, since $a||b$ in~$\P$ we now know that $a >_{\overline{\sigma}} b$. Because $a = \max_{\overline{\sigma}} A$ and $\BD \setminus \AD =\{b\}$, we have $ \max_{\overline{\sigma}} (\AD \triangle \BD)=a$.  

Let us now consider the class $\mathcal{C}_{D,I}$ defined by $D = A \triangle B$ and $I =A \cap B$. The elements $a$ and $b$ lie in different components of~$\P[D]$, because $B \setminus A = \{b\}$ and $a ||b$ in $\P$. So we may assume that $a \in \alpha$ and $b \in \beta$,  where $\alpha$ and $\beta$ are different elements from the set $\mc{K}$ of components of $\P[D]$. As before, set $d=|\mc{K}|$.

Observe that $|\AD \cup \BD|=|\AD \cup \{b\}|=k+1$. Now choose a subset \mbox{$K \subset \mc{K}$} with $\alpha, \beta \notin K$. For the corresponding downset $(X_K \cup I)^{\da}= A_K^{\da}$ we have $ A_K^{\da} \subseteq \AD \cup \BD$. Since $a,b \notin A_K^{\da}$, we can apply the induction hypothesis to $A_K^{\da}$. We can also apply it to the set $(b \cup I)^{\da}$. It holds that $\max_\sigma (A_K \triangle (b\cup I))=b$, so we have $ A_K^{\da} < (b \cup I)^{\da}$ in $L$ by induction. Let $X_K^c \cup I = B_K$, then we have $(b \cup I)^{\da} \subseteq B_K^{\da}$, and thus $A_K^{\da} < B_K^{\da}$ in $L$ by transitivity. Analogously,  $\max_{\overline{\sigma}} (A_K \triangle (a\cup I)=a$ holds, so we have $A_K^{\da} < (a \cup I)^{\downarrow}$ in $\overline{L}$ by induction and thus $A_K^{\da}  <(a \cup I)^{\da} < (X_K^c \cup I)^{\da} = B_K^{\da}$ in $\overline{L}$. 

It follows that $\AD_K$ is down in $L$ and~$\overline{L}$ for every $K \subset \mc{K}$ with $\alpha, \beta \notin K$.
Thus $(\AD_K, \BD_K)$ cannot yield a reversal between $L$ and~$\overline{L}$, and neither can $(\BD_K, \AD_K)$. There are $2^{d-2}$ possibilities to choose $K$. Thus we have exhibited $2 \cdot 2^{d-2}$ pairs in $\mathcal{C}_{D,I}$ which do not contribute a reversal between $L$ and $\overline{L}$. From the fact we remarked after Theorem~\ref{thm_Lsigma_diam} it follows that all other pairs in $\mathcal{C}_{D,I}$ have to contribute reversals. 

Consequently, all subsets $K$ of $\mc{K}$ containing exactly one of the two components $\alpha$ and $\beta$ need to contribute a reversal, or equivalently, all $A_K^{\da}$ which contain exactly one of the two elements $a,b$ need to be down in exactly one of the two linear extensions. In particular, our set $\AD$ needs to be down relative to~$\C_{D,I}$ in exactly one of the two linear extensions. 

 It turns out that $\AD$ cannot be down in $\overline{L}$: For \mbox{$(I \cup a)^{\downarrow} \subset \AD$}, we have $\max_{\overline{\sigma}} (B \triangle (I \cup a)) = a $, so we have  $\BD < (I \cup a)^{\da}$ in $\overline{L}$ by induction and thus $\BD< \AD$ in $\overline{L}$ by transitivity. Hence it follows that $\AD < \BD$ in $L$. This proves our claim that $\AD$ has to sit in the slot between $\tA$ and $\BD$ in $L^k$. 
\qedclaim

Because $\max_{\sigma} (A \triangle B)=b \in B$, and $\AD < \BD$ in $L$ as shown in the claim, we now know that $\AD$ is in $\sigma$-revlex order with $\BD$ in $L^k$. Since the slot after~$\tA$ is the lowest possible position for $\AD$ in $L^k$, it follows from the transitivity of the $\sigma$-revlex order that~$\AD$ is in $\sigma$-revlex order in $L^k$ with all sets of smaller cardinality. By reversing the roles of $L$ and $\overline{L}$, we obtain that~$\AD$ is in $\overline{\sigma}$-revlex order in $\overline{L}^k$ with all sets of smaller cardinality.  Next we show that all pairs of sets with equal cardinality $k$ are in $\sigma$-revlex order in~$L^k$. 

Let $A_i^{\da},\, A_j^{\da} \in \mathcal{D}_\P$ be two downsets of the same cardinality $k$. If they are inserted into different slots in $L$, they are in $\sigma$-revlex order by transitivity.  If they belong into the same slot, this means that they have the same largest $(k-1)$-subset $\tA$ in $L$.  So their symmetric difference contains only two elements, say, $A_i^{\da} \triangle A_j^{\da} = \{a_i, a_j\}$. We have \mbox{$a_i= A_i^{\da} \setminus \tA = \min_\sigma A_i = \max_{\overline{\sigma}} A_i $} and \mbox{$a_j=A_j^{\da} \setminus \tA =\min_\sigma A_j = \max_{\overline{\sigma}} A_j$}. Note that $a_i$ and $a_j$ have to be incomparable in~$\P$, and assume that $a_i<_\sigma a_j$, thus $a_j <_{\overline{\sigma}} a_i$. Then for the pair $\{A_j^{\da}, a_i^{\da}\}$ we know $\max_{\overline{\sigma}} ( A_j^{\da} \triangle a_i^{\da})=a_i$. Hence by induction, $A_j^{\da} < a_i^{\da} < A_i^{\da}$ in~$\overline{L}$. But since the equivalence class containing $(A_i^{\da}, A_j^{\da})$ needs to contribute $2^{2-2}=1$ reversal between $L$ and $\overline{L}$, and this can only be the pair $\{A_i^{\da}, A_j^{\da}\}$, we know that we must have $A_i^{\da} < A_j^{\da}$ in $L$. Because of $\max_\sigma A_i \triangle A_j = a_j$, this means that $A_i^{\da}$ and $A_j^{\da}$ are in revlex order in $L$.

We can apply the same argument (only with the roles of $L$ and $\overline{L}$ reversed) to show that all pairs of downsets of equal cardinality $k$ are in $\overline{\sigma}$-revlex order in $\overline{L}$. By induction we conclude that $L=L_\sigma$ and $\overline{L} = L_{\overline{\sigma}}$. \qed

\subsection{Computing the Linear Extension Diameter}
\label{subsec_latdown_comp}

It is NP-complete to compute the linear extension diameter of a
general poset, see~\cite{massow-diam}. That is, the linear
extension diameter of a general poset~$\P$ cannot be
computed in time polynomial in~$|\P|$ (unless $\textmd{P}=\textmd{NP}$). But with the results of the
previous section, the problem is tractable if the poset is a downset
lattice~$\mc{D}$ of a 2-dimensional poset. 

In fact, we can construct any diametral pair of linear extensions of $\mc{D}$ in time
polynomial in $|\mc{D}|$. To see this, we use a well-known fact about
distributive lattices: From a downset lattice~$\mathcal{D}$ one can obtain~$\P$ with $\mathcal{D}=\mathcal{D}_\P$ as the poset induced by the join irreducible elements of~$\mathcal{D}$. 
Finding a minimal realizer $\sigma, \overline{\sigma}$ of the 2-dimensional poset $\P$ amounts to finding a transitive orientation of its incomparability graph (cf.~e.g.~\cite{moehring-tractable}). This can be done in time linear in~$|\P|$~\cite{mcconnell-spinrad}. Then $\sigma, \overline{\sigma}$ is a diametral pair  of linear extensions of~$\P$. With the definition of the $\sigma$-revlex order we can compute~$L_\sigma$ and~$L_{\overline{\sigma}}$, and this is a diametral pair of linear extensions of~$\mathcal{D}$ by Theorem~\ref{thm_Lsigma_diam}. We know by Theorem~\ref{thm_unique_downlat} that all diametral pairs arise this way. The linear extension diameter of $\mc{D}$ can now be computed by just checking for all pairs of elements of $\mc{D}$ whether they form a reversal between $L_\sigma$ and~$L_{\overline{\sigma}}$.

In this subsection we show that we can in fact do much better: For a \mbox{2-dimensional} poset~$\P$, we
can compute the linear extension diameter of $\mc{D}_\P$ in time
polynomial in $|\P|$. Note that in general,
$\P$ can have exponentially many downsets, so $|\mc{D}_\P|$ is
exponentially larger than $|\P|$. 

For the proofs of this subsection, we mainly consider antichains instead 
of downsets, again using the canonical bijection between them. 
It is known that the antichains of a 2-dimensional poset can be
counted in polynomial time, see~\cite{steiner} or~\cite{moehring-tractable}. We
give a proof in the lemma below since the methods we use for the following
theorem rely on the same ideas.

\begin{lemma}
\label{lem_count_antichains}
Let $\P$ be a 2-dimensional poset. Denote by $A(\P)$ the set of
antichains of $\P$ and let $a(\P)=|A(\P)|$. Then $a(\P)$ can be
computed in time~$\mc{O}(|\P|^2)$.
\end{lemma}

\proof Let $\sigma=x_1 x_2 \ldots x_n$ be a non-separating linear
extension of $\P$. Denote by $A(x_i)$ the set of antichains of $\P$
which contain $x_i$ as $\sigma$-largest element, and let
$a(x_i)=|A(x_i)|$. We will use a dynamic programming approach to
compute $a(x_i)$ for all $i$.

To start with, we have $a(x_1)=1$. Now suppose we have computed $a(x_j)$
for all $j<i$. The main observation is that for any $A \in A(x_j)$
with $j<i$ and $x_i || x_j$, the set $x_i \cup A$ is again an
antichain. This holds because any $x_k \in A$ with $x_k < x_i$ would
yield a contradiction to $\sigma$ being non-separating. Therefore we
have
\[
a(x_i)= 1+ \sum_{j<i, \; x_i || x_j} a(x_j),
\]
where the 1 accounts for the antichain $\{x_i\}$. Consequently, we obtain
the number of all antichains of $\P$ as $a(\P)= 1 + \sum_i a(x_i)$, where 
the 1 accounts for the empty set.

With the above formula, the evaluation of $a(x_i)$ can be done in 
linear time for each $i$.
Thus $a(\P)$ can be computed in quadratic time.  \qed

\begin{theorem}
The linear extension diameter of the downset lattice $\mc{D}_\P$ of
a 2-dimensional poset $\P$ can be computed in time $\mc{O}(|\P|^5)$.
\end{theorem}

\proof From Corollary~\ref{corr_led_downlat}, we know that
$\textmd{led}(\mc{D}_\P)$ equals a quarter of the number of pairs
$(A, B)$ of antichains of $\P$ such that $\P[A \triangle B]$ has at least two
connected components. For a pair $(A,B)$ of antichains of $\P$, 
we set $D=A \triangle B$ and $I = A \cap B$. 

We will count four different classes of pairs of antichains. Let
$\alpha$ be the number of \emph{all} ordered pairs of
antichains of $\P$. Let $\beta$ be the number of pairs $(A, B)$
with $D = \emptyset$, and $\gamma$ the number of pairs with $|D|=1$.
Finally, let $\delta$ be the number of pairs such that $|D|>1$ and
$\P[D]$ is connected. Then $\textmd{led}(\mc{D}_\P)=
\frac{1}{4}\left(\alpha - \beta - \gamma - \delta\right)$.

We have $\alpha= a(\P)^2$. The pairs we count for $\beta$ are just
the pairs $(A, A)$, so $\beta= a(\P)$. 

For the following,  let $\sigma=x_1 x_2 \ldots x_n$ be a non-separating 
linear extension of $\P$. We denote by $[x_i, x_k]$ the set $\{x_i, x_{i+1}, \ldots,
x_k\}$, and by $(x_i, x_k)$ the set $\{x_{i+1}, \ldots, x_{k-1}\}$. We 
use analogous notions for ``half-open intervals'' of~$\sigma$.

To obtain $\gamma$, we count the pairs $(A, A-x)$, where $A$ is a
non-empty antichain in $\P$, and $x$ is an element of $A$. Thus we
want to count each $A$ exactly $|A|$ times. We want to refine the
ideas of the proof of Lemma~\ref{lem_count_antichains} to keep track
of the sizes of the antichains. Hence we define vectors $s(x_i)$,
where~$s_r(x_i)$ is the number of antichains of cardinality $r$ in
$A(x_i)$.

We can recursively compute $s(x_i)$ for \mbox{$i=1,2, \ldots, n$} as
follows: The first entry of each $s(x_i)$ is 1, counting the
antichain $\{x_i\}$. For the other entries we have
\[
s_r(x_i)= \sum_{j<i: \, x_j ||x_i} s_{r-1}(x_j).
\]
Then the number of pairs $(A, A-x)$ equals $\sum_{r} r \sum_i
s_r(x_i)$. Now $\gamma$ is twice this number, since we also need to
count the pairs $(A-x, A)$. 

\begin{figure}[ht]
  \begin{center}
    \psfrag{D}{$D$}
    \psfrag{Il}{$I^{\txt{left}}$}
    \psfrag{Ir}{$I^{\txt{right}}$}
    \includegraphics[scale=0.5]{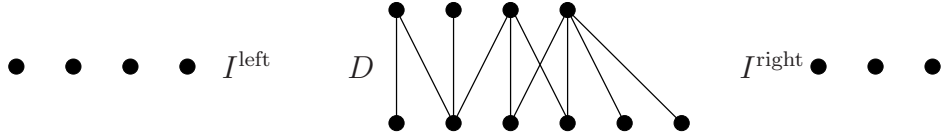}
    \label{fig_delta}
\caption{$\P[D \cup I]$ for a pair of antichains counted in $\delta$.}
   \end{center}
\end{figure}

The most difficult part is to compute $\delta$, the number of pairs $(A,B)$ 
of antichains of $\P$ such that $|D|>1$ and $\P[D]$ is connected. 
Let us take a look at the structure of $\P[D\cup I]$, see~Figure~\ref{fig_delta}. 

We know that~$\P[I]$ is an antichain and $\P[D]$ consists of two antichains, $\Max(D)$ and $\Min(D)$, 
which are both non-empty by definition of $\delta$.  
We claim that each element of $I$ is either $\sigma$-smaller than all elements of $D$, 
or $\sigma$-larger than all elements of $D$. Indeed, suppose there is an $x \in I$ and 
$u,v \in D$ with $u<_\sigma x<_\sigma v$. Then since $\P[D]$ is connected, we can find 
$u',v' \in D$ with $u' < v'$ in $\P$ and $u'<_\sigma x<_\sigma v'$. But since there are no relations 
between $x$ and the elements of $D$, this means that $\sigma$ is 
separating. This contradiction proves our claim. 
Thus $\sigma$ can be partitioned into three intervals, such that the elements of~$D$ are 
all contained in the middle interval, and $I$ is split up into two parts: $I^{\textmd{left}}$, 
living on the first interval, and $I^{\textmd{right}}$, living on the third interval. 


Now define $\delta(k,\ell)$ as the number of pairs counted for $\delta$ which fulfill 
\mbox{$\max_\sigma \Min(D)=x_k$} and \mbox{$\max_\sigma \Max(D) =
x_\ell$}. In addition, we require that \mbox{$x_\ell = \max_\sigma \A \cup \B$}, which means that $I^{\textmd{right}}$ is empty. 
We split up $\delta(k,\ell)$ into the number $\delta_1(k, \ell)$ of 
pairs for which $\P[D]$ has only one maximum and the number 
$\delta_2(k, \ell)$ of pairs for which it has several. 

\begin{figure}[h]
  \begin{center}
    \psfrag{x1}{$x_1$}
    \psfrag{x2}{$x_2$}
    \psfrag{x3}{$x_3$}
    \psfrag{xk}{$x_k$}
    \psfrag{xl}{$x_\ell$}
    \psfrag{xi}{$x_i$}
    \psfrag{D}{$D$}
    \psfrag{I}{$I^{\txt{left}}$}
    \includegraphics[scale=0.5]{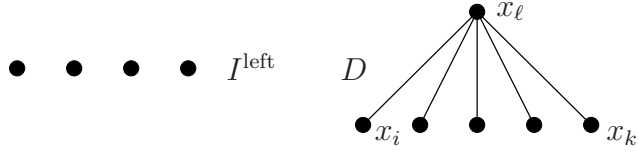}
    \label{fig_deltakl1}
\caption{$\P[D \cup I]$ for a pair of antichains counted in $\delta_1(k,\ell)$.}
   \end{center}
\end{figure}

To compute $\delta_1(k, \ell)$, we have to count the number of possibilities to 
choose the antichain~$\Min(D)$ and the antichain~$I^{\txt{left}}$. By definition we have \mbox{$\max_\sigma
\Min(D)=x_k$}. Suppose that $\min_\sigma \Min(D) = x_i$ as in ~Figure~\ref{fig_deltakl1}. 

Let us define $\P_{i,k,\ell}$ as the poset induced by the elements $x_j \in (x_i, x_k)$ with $x_j <
x_\ell$ and $\{x_i, x_j, x_k\} \in A(\P)$. Then to choose $\Min(D)$,
we have to choose an antichain in $\P_{i,k,\ell}$.

Once a set $D$ is fixed, it remains to choose $I^{\txt{left}}$ to determine a pair of 
antichains counted in $\delta_1(k, \ell)$. Let $\P_{i, \ell}^{\textmd{left}}$ be the poset induced by the
elements $x \in \P$ with $x \in [x_1, x_i)$ and $x || x_\ell$. We claim that the sets which can be chosen 
as $I^{\txt{left}}$ are exactly the antichains of $\P_{i, \ell}^{\textmd{left}}$.

By definition, each element $x \in I^{\txt{left}}$ is $\sigma$-smaller than $x_i$. 
We have to choose $x$ so that it is incomparable to all elements in $D$. 
It is clear that $x \in [x_1, x_i)$ cannot be larger in $\P$ than any element
in $D \subseteq [x_i, x_\ell]$. Now if we choose $x$ incomparable to $x_\ell$, it cannot be smaller than any
 element in~$D$, either. Thus to choose $I$, we have to choose an antichain in 
$\P_{i,\ell}^{\textmd{left}}$ as claimed. 
Altogether we have
\[
    \delta_1(k, \ell) = \sum_{x_i \in [x_1, x_k], \; \{x_i, x_k \} \in A(\P), \;
    x_i < x_\ell} a(\P_{i, k, \ell}) \cdot a(\P_{i, \ell}^{\textmd{left}}).
\]

\begin{figure}[h]
  \begin{center}
    \psfrag{x1}{$x_1$}
    \psfrag{x2}{$x_2$}
    \psfrag{x3}{$x_3$}
    \psfrag{xk}{$x_k$}
    \psfrag{xl}{$x_\ell$}
    \psfrag{xk'}{$x_{k'}$}
    \psfrag{xl'}{$x_{\ell'}$}
    \psfrag{D}{$D$}
    \psfrag{I}{$I^{\txt{left}}$}
    \includegraphics[scale=0.5]{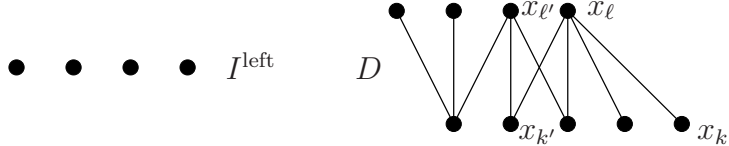}
\caption{$\P[D\cup I]$ for a pair of antichains counted in $\delta_2(k, \ell)$.}
    \label{fig_deltakl2}
   \end{center}
\end{figure}

It remains to compute $\delta_2(k, \ell)$, the number of pairs $(A,B)$ counted in 
$\delta(k, \ell)$ for which $\P[D]$ has several maxima. We want to
cut off the $\sigma$-largest maximum and (possibly) some minima and recursively
use values $\delta_2(k',\ell')$ and $\delta_1(k',\ell')$ that we have calculated already 
(cf.~Figure~\ref{fig_deltakl2}). 

For a pair $(A, B)$ counted in $\delta_2(k, \ell)$, the second largest maximum of~$\P[D]$ in 
$\sigma$ is an element $x_{\ell'} \in [x_1,x_\ell)$ with $x_{\ell'} || x_\ell $. In general,
$\P[D]$ is not connected after deletion of $x_\ell$. But since $\P[D]$ was connected originally, $\Min(D)$ contains an element~$x_{k'}$ with $x_{k'} < x_\ell$ and 
$x_{k'} < x_{\ell'}$. Let $x_{k'}$ be the $\sigma$-smallest such element. 
There can be more elements in $[x_{k'}, x_k]$ which are part of $\Min(D)$. With the same reasoning
as for  $\delta_1(k, \ell)$, these form exactly the antichains in~$\P_{k', k, \ell}$. 
So we have
\[
    \delta_2(k, \ell)  = \sum_{x_{\ell'} \in [x_1, x_\ell), \; x_{\ell'}|| x_\ell} \;\,
\sum_{x_{k'} \in [x_1, x_{\ell'}), \; x_{k'}<x_{\ell'}, x_\ell}
 \big( \delta_1(k', \ell')+\delta_2(k', \ell') \big) 
\cdot a(\P_{k', k, \ell}).
\]

Recall that for the pairs counted in $\delta(k,\ell)$, we required that 
$I^{\txt{right}}$ is empty. So to compute $\delta$, we have to weight every pair with the number of possible choices for $I^{\txt{right}}$. Let $\P_{k,\ell}^{\txt{right}}$ be the poset induced by the elements of~$\P$ which are in $(x_\ell, x_n]$ and incomparable to $x_k$ in $\P$. We claim that the sets eligible for $I^{\txt{right}}$ are exactly the antichains of~$\P_{k,\ell}^{\txt{right}}$. 

Each element $x \in I^{\txt{right}}$ has to be incomparable to all elements of $D$. 
Since $x \in (x_\ell, x_n]$, it cannot be smaller than any element in $D$. If we choose~$x$
incomparable to $x_k$, it cannot be larger than any element of $D$ either: If $x > y$ for 
some element $y \in \Min(D)$, then $y$ has to be $\sigma$-smaller than 
$x_k$, which makes $\sigma$ separating. Thus to choose $I^{\txt{right}}$ we have to choose 
an antichain in~$\P_{k,\ell}^{\textmd{right}}$ as claimed. Hence we obtain $\delta$ as follows:
\[
    \delta = \sum_{k, \ell}      \big( \delta_1(k, \ell) + \delta_2(k, \ell) \big) \cdot a(\P_{k, \ell}^{\textmd{right}}).
\]

What is the overall running time for the computation of $\txt{led}(\mc{D}_\P)$? From Lemma~\ref{lem_count_antichains} we know that the number of antichains of a poset can be 
computed in quadratic time. Thus $\alpha$ and $\beta$ can be determined in quadratic time.
To compute $\gamma$, we need to compute $s_r(x_i)$ for $r=1\ldots n$ and $i=1 \ldots n$. 
For each value $s_r(x_i)$, our formula can be evaluated in linear time. Thus it takes $O(n^3)$ 
to determine $\gamma$. 

For the computation of $\delta$ we first determine in a preprocessing step the \mbox{values}~$a(\P_{i, k, \ell})$ for all triples~$i,k, \ell$. Given such a triple, we can build~$\P_{i, k, \ell}$ in linear 
time, and then compute~$a(\P_{i, k, \ell})$ in quadratic time using Lemma~\ref{lem_count_antichains} again. 
Altogether this can be done in~$O(n^5)$. Similarly, we can \mbox{determine} $a(\P_{k, \ell}^{\textmd{left}})$ and 
$a(\P_{k, \ell}^{\textmd{right}})$ for all pairs~$k, \ell$ in a preprocessing step taking time~$O(n^4)$. 
Then for each pair~$k,\ell$, we can compute~$\delta_1(k, \ell)$ and~$\delta_2(k, \ell)$ in linear time. 
Thus it takes~$O(n^3)$ to obtain all these values. In the end, we can put them together in quadratic time to obtain~$\delta$. 

The overall running time is the maximum over all these separate steps. We conclude that $\txt{led}(D_\P)$ can be computed in time $O(n^5)$. 
\qed

In the previous theorem we showed how to compute $\textmd{led}(D_\P)$ for \mbox{2-dimensional $\P$}, but we could not give an explicit formula like the one we have for the Boolean lattice. This is possible for the special case where 
$\P$ is a disjoint union of chains. These lattices $\mathcal{D}_\P$ are also known as the factor lattices of integers: If $\P= C_1 \cup \ldots \cup C_w$ with $|C_i|=\ell_i$, we can associate each chain with a prime number $p_i$, and then $\mathcal{D}_\P$ is the lattice of all factors of $m=\prod_{i=1}^w p_i^{\ell_i}$, ordered by divisibility. 

\begin{theorem}
\label{thm_led_factorlattice}
	If $\P =C_1 \cup \ldots \cup C_w$ with $|C_i|=\ell_i$ is a disjoint union of chains, then the linear extension diameter of $\mathcal{D}_\P$ equals
\[
	\frac{1}{4} \cdot \left( \left( \prod_{i=1}^{\omega} (\ell_i +1) \right)^2 - \sum_{k=1}^{\omega} (\ell_k +1)\ell_k \cdot \prod_{i \neq k} (\ell_i +1) - \prod_{i=1}^{\omega} (\ell_i +1)\right).
\]
\end{theorem}

\proof
From Corollary~\ref{corr_led_downlat} we know that $\textmd{led}(\mathcal{D}_\P)$ equals a quarter of the number of pairs $(A, B)$ of antichains of $\P$ such that $\P[A \triangle B]$ has at least two connected components. So we need to count the pairs of antichains of~$\P$ which differ on at least two of the $C_i$. We will count \emph{all} pairs of antichains and substract from it the number of pairs differing on zero or one chain. 

To choose one antichain, we have $\ell_i +1$ choices in each $C_i$. So $\P$ contains exactly $ \prod_{i=1}^{\omega} (\ell_i +1)$ antichains, and this is also the number of pairs of antichains differing on zero chains. The number of all pairs of antichains is thus $\left( \prod_{i=1}^{\omega} (\ell_i +1) \right)^2$. The number of pairs of antichains which differ on one chain is the sum over $k$ of all choices of two different elements in chain $C_k$ and one element from each other chain. This yields the desired formula. \qed


\section{Open Problems}

It is NP-hard to compute the linear extension diameter of a general poset, see~\cite{massow-diam}. For Boolean lattices and for downset lattices of 2-dimensional posets we can now construct the diametral pairs of linear extensions in polynomial time. Is this possible for arbitrary distributive lattices? 

\begin{open}
	Is it possible to compute the linear extension diameter of an arbitrary distributive lattice in polynomial time, or even characterize its diametral pairs of linear extensions?
\end{open}

Another question we are interested in asks whether there is a fixed fraction of the incomparable elements of a poset that can always be reversed between two linear extensions. 

\begin{open}
	Is there a constant $c>0$ such that $\textmd{led}(\P)/ \textmd{inc}(\P)>c$ for all posets $\P$?
\end{open}

It is also interesting to look at subposets of the Boolean lattice. Consider for example the subposet of $B_5$ induced by the sets of cardinality 2 and 3. It turns out that the revlex linear extensions do not form a diametral pair of linear extensions  of this poset. 

\begin{open}
 Can we construct the diametral pairs of linear extensions of a poset induced by two levels of~$B_n$?
\end{open}

\noindent {\em Acknowledgement.} We thank Kolja Knauer for fruitful discussions. 

\bibliography{my_library}
\bibliographystyle{plain}

\end{document}